\documentclass[doublespacing]{elsart3-1}

\usepackage{amssymb}
\usepackage{calc,amsfonts,amscd,epsfig,psfrag,amsmath,amssymb,enumerate,graphicx,mathrsfs,paralist}
\usepackage[english,francais]{babel}

\newtheorem{theorem}{Theorem}[section]
\newtheorem{maintheorem}{Theorem}
\newtheorem{maincorollary}[maintheorem]{Corollary}

\newtheorem{e-proposition}[theorem]{Proposition}

\newtheorem{e-definition}[theorem]{Definition\rm}

\newtheorem{question}[theorem]{Question}

\renewcommand{\theequation}{\arabic{equation}}
\newcommand{\Dim}{\textsc{d} }
\newcommand{\Z}{\mathbb Z}

\newcommand{\sL}{\ensuremath{\mathscr L}}
\newcommand{\sW}{\ensuremath{\mathscr W}}

\renewcommand{\epsilon}{\varepsilon}

\def\og{\leavevmode\raise.3ex\hbox{$\scriptscriptstyle\langle\!\langle$~}}
\def\fg{\leavevmode\raise.3ex\hbox{~$\!\scriptscriptstyle\,\rangle\!\rangle$}}

\journal{the Acad\'emie des sciences}
\begin{document}
\centerline{}
\begin{frontmatter}


\selectlanguage{english}
\title{The shape of the $(2+1)$D SOS surface above a wall}


\selectlanguage{english}

\author[authorlabel1]{Pietro Caputo}~~
\ead{caputo@mat.uniroma3.it}
\author[authorlabel2]{Eyal Lubetzky}~~
\ead{eyal@microsoft.com}
\author[authorlabel1]{Fabio Martinelli}~~
\ead{martin@mat.uniroma3.it}
\author[authorlabel3]{Allan Sly}~~
\ead{sly@stat.berkeley.edu}
\author[authorlabel4]{Fabio Lucio Toninelli}
\ead{fabio-lucio.toninelli@ens-lyon.fr}
\thanks[thankslabel1]{This work was supported by the European Research Council through the ``Advanced
Grant'' PTRELSS 228032}
\address[authorlabel1]{Universit\`a Roma Tre, Largo S.\ Murialdo 1, 00146 Roma, Italia.}
\address[authorlabel2]{Microsoft Research, One Microsoft Way, Redmond, WA 98052, USA.}
\address[authorlabel3]{UC Berkeley, Berkeley, CA 94720, USA.}
\address[authorlabel4]{CNRS and ENS Lyon, Laboratoire de Physique, 46 All\'ee d'Italie, 69364 Lyon, France.}


\medskip
\begin{center}
\end{center}

\begin{abstract}
\selectlanguage{english}
We give a full description for the shape of the classical (2+1)\Dim Solid-On-Solid model above a wall, introduced by
Temperley (1952). On an $L\times L$ box at a large inverse-temperature $\beta$ the height of most sites concentrates on a single level $h = \lfloor \frac1{4\beta}\log L\rfloor$ for most values of $L$.
For a sequence of diverging boxes the ensemble of level lines of heights $(h,h-1,\ldots)$ has a scaling limit in Hausdorff distance iff the fractional parts of $\frac1{4\beta}\log L$ converge to a noncritical value. The scaling limit is explicitly given by  nested distinct loops formed
via translates of Wulff shapes.
 Finally, the $h$-level lines feature $L^{1/3+o(1)}$ fluctuations from the side boundaries.

\vskip 0.5\baselineskip

\selectlanguage{francais}
\noindent{\bf R\'esum\'e} \vskip 0.5\baselineskip \noindent
{\bf La forme de l'interface SOS (2+1)-dimensionnelle au-dessus d'un mur.}
Nous donnons une description compl\`ete de la forme typique de l'interface SOS (2+1)-dimensionnelle au-dessus d'un mur, introduite par Temperley (1952). Dans une bo\^ite $L\times L$ \`a basse temp\'erature  $T=1/\beta$, la hauteur de la plupart des sommets se concentre au niveau $h =  \lfloor \frac1{4\beta}\log L\rfloor$, pour la plupart des valeurs de $L$. Pour une suite croissante de bo\^ites, l'ensemble de lignes de niveau de hauteur $(h,h-1,\ldots)$ admet une limite au sense de la distance de Hausdorff ssi la partie fractionnaire de  $\frac1{4\beta}\log L$ converge \`a une valeur non critique.
La limite d'\'echelle est donn\'ee explicitement par des boucles imbriqu\'ees, form\'ees par des translations de formes de Wulff.
Enfin, la distance entre le bord de la bo\^ite et la ligne de niveau $h$ a des fluctuations $L^{1/3+o(1)}$.

\end{abstract}
\end{frontmatter}


\selectlanguage{english}
\section{Introduction}

The $(2+1)$-dimensional \emph{Solid-On-Solid} model is a crystal surface model whose definition goes back to Temperley~\cite{Temperley} in 1952 (also known as the Onsager-Temperley sheet).
At low temperatures, the model approximates the interface between the plus and minus phases in the $3$\Dim Ising model.
See, e.g.,~\cite{FS1,FS2,FS3} for notable previous works on this model.

The configuration space of the model on a finite box $\Lambda\subset\Z^2$ with a floor (wall) at $0$ and zero boundary conditions is the set $\Omega_{\Lambda}$ of all height functions $\eta$ on $\Z^2$ such that $\Lambda\ni x \mapsto \eta_x \in \Z_+$ whereas $\eta_x=0$ for all $x\notin\Lambda$.
The probability of $\eta\in\Omega_{\Lambda}$ is given by the Gibbs distribution
\begin{equation}
  \label{eq-ASOS}
  \pi_\Lambda(\eta) = \frac{1}{Z_\Lambda} \exp\bigg(-\beta \sum_{x\sim y}|\eta_x-\eta_y|\bigg)\,,
\end{equation}
where $\beta>0$ is the inverse-temperature, $x\sim y$ denotes a nearest-neighbor bond in the lattice $\Z^2$ and the normalizing constant $Z_\Lambda$ is the partition-function.

Despite progress on understanding the typical height of the surface
(notably the work of Bricmont, El-Mellouki and Fr{\"o}hlich~\cite{BEF} in 1986 and the recent companion paper of the authors~\cite{CLMST}),
little was known on its actual 3\Dim shape. The fundamental problem is the following:
\begin{question}
Consider the ensemble of all level lines of the low temperature $(2+1)$\Dim SOS on an $L\times L$ box, rescaled to the unit square.
\begin{compactenum}[(i)]
\item \label{q-item-1} Do these jointly converge to a scaling limit as $L\to\infty$, e.g., in Hausdorff distance?
    \item \label{q-item-2} If so, can the limit be explicitly described?
    \item \label{q-item-3} For finite large $L$, what are the fluctuations of the level lines around their limit?
\end{compactenum}
\end{question}
\medskip
In this work we fully resolve parts~\eqref{q-item-1} and~\eqref{q-item-2} and partially answer part~\eqref{q-item-3}. En route, we also establish that for most values of $L$ the surface height concentrates on the single level $\lfloor \frac1{4\beta}\log L\rfloor$.

\subsubsection*{Main results}
Our first result addresses the distribution of the surface height.

\medskip
\begin{maintheorem}[Height Concentration]\label{mainthm-1}
Fix $\beta > 0$ to be sufficiently large and define
\begin{equation}\label{eq-H-def}
H(L) = (1/4\beta) \log L\,.
\end{equation}
Let $E_h = \left\{ \eta : \#\{x : \eta_x=h\} \geq \frac9{10}L^2\right\}$ be the event that at least $\frac9{10}$ of the sites are at height $h$.
Then the SOS measure on the $L\times L$ box $\Lambda\subset \Z^2$ with inverse-temperature $\beta$ satisfies
\begin{equation}
  \label{eq-two-levels}
  \lim_{L\to\infty} \pi_\Lambda\left(E_{\lfloor H \rfloor - 1} \cup E_{\lfloor H\rfloor}\right) = 1\,.
\end{equation}
Furthermore, the typical height of the configuration is governed by $\alpha(L) = H - \lfloor H \rfloor$
as follows:
  Let $\Lambda_n$ be a diverging sequence of boxes with side-lengths $L_n$. There exists $0<\alpha_c(\beta) <1$ so that
\begin{enumerate}
  [(i)]
  \item If $\liminf_{n\to\infty}\alpha(L_n) > \alpha_c$ then $\lim_{n\to\infty} \pi_{\Lambda_n}\left(E_{\lfloor H\rfloor}\right) = 1$.
  \item If $\limsup_{n\to\infty}\alpha(L_n) < \alpha_c$ then $\lim_{n\to\infty} \pi_{\Lambda_n}\left(E_{\lfloor H\rfloor-1}\right) = 1$.
\end{enumerate}
\end{maintheorem}
\medskip
\noindent {\bf Remark.}
The constant $\frac{9}{10}$ in the definition of $E_h$ can be replaced by $1-\epsilon$ for any arbitrarily small $\epsilon>0$ provided that $\beta$ is chosen large enough.
Moreover, the critical fractional value $\alpha_c(\beta)$ satisfies $\alpha_c = (1+\epsilon_\beta) \frac{\log(4\beta)}{4\beta} $ where $\epsilon_\beta \to 0$ as $\beta\to\infty$. Therefore, for large enough (fixed) $\beta$ most values of $\alpha(L)$ will yield $\pi_\Lambda(E_{\lfloor H \rfloor})\to 1$.

We now address the scaling limit of the ensemble of level lines, formally defined next.

\begin{e-definition}
Let $\eta$ be an SOS configuration on a box $\Lambda\in\Z^2$. The $h$-level lines ($h\geq 1$) are a collection of disjoint self-avoiding loops in the dual $\Z^{2*}$ formed as follows.
Let $\mathcal{E}_h$ be the set of all bonds $e' \in \Z^{2*}$ whose dual edge $e=(x,y)$ is such
that $\eta_x \geq h$ whereas $\eta_y < h$, or vice versa.
Due to the zero boundary conditions on $\eta$, the bonds of $\mathcal{E}_h$ can then be uniquely partitioned into a finite number of edge-disjoint simple loops as follows. Whenever four bonds in $\mathcal{E}_h$ meet at a vertex we separate them along the NW-oriented diagonal going through the intersection. The final collection of loops comprises the $h$-level lines.
\end{e-definition}

The next theorem gives a necessary and sufficient condition for the existence of a scaling limit to the full ensemble of level lines, and furthermore provides an explicit description for the limit.

\medskip
\begin{maintheorem}[Shape Theorem]\label{mainthm-2}
 Fix $\beta > 0$ to be sufficiently large and let $L_n$ be a diverging sequence of
 side-lengths. Set $H_n=(1/4\beta)\log L_n$ and $\alpha_n = H_n - \lfloor H_n\rfloor$. For an SOS surface on a box $\Lambda_n$ of side-length $L_n$, let $(\sL_0^{(n)},\sL_1^{(n)},\ldots)$ be the collections of loops whose area is at least $\log^2 L$
belonging to the level lines of heights $(\lfloor H_n \rfloor,\lfloor H_n \rfloor -1,\ldots)$, respectively. Then:

\begin{compactenum}
  [(a)]
\item With high probability, the level lines of every height $h>\lfloor H \rfloor$ consist
of microscopic loops, each of which spans an area of at most $\log^2 L$.
In addition, for each $i\geq 1$ w.h.p.\ $\sL_i^{(n)}$ consists of exactly one loop whereas $\sL_0^{(n)}$ is w.h.p.\ either empty or contains a single loop.

\item The rescaled loop ensemble $\frac{1}{L_n}(\sL_0^{(n)},\sL_1^{(n)},\ldots)$ converges almost surely to a limit $(\sW_0,\sW_1,\ldots)$ in Hausdorff distance iff $\lim_{n\to\infty} \alpha_n$ exists and differs from $\alpha_c$, the critical fractional value given in Theorem~\ref{mainthm-1}. More precisely, in this case almost surely
    $ \limsup_{n\to\infty}\sup_i d_\mathcal{H}(\frac{1}{L_n}\sL_i^{(n)},\sW_i) = 0$.
\item The scaling limit $(\sW_0,\sW_1,\ldots)$ is explicitly given as a nested sequence of distinct loops invariant under $\pi/4$ rotations, each of which is constructed as follows. Let $\alpha_\star = \lim_{n\to\infty}\alpha_n$.

    \begin{compactitem}
      \item If $\alpha_\star > \alpha_c$ then for any $i\geq 0$ the loop $\sW_i$
    is the boundary of the union of all possible translates within the unit square
    of $\sW_\star(r_i)$ with an explicit $0<r_i(\alpha_\star) <1$,
    where $\sW_\star(r)$ is the $r$-dilation of a convex smooth shape (the area 1 SOS Wulff shape).

\item
    If $\alpha_\star < \alpha_c$ then $\sW_0$ is empty and the remaining loops $\{\sW_i\}_{i\geq 1}$ are constructed as above.
    \end{compactitem}
\end{compactenum}
\end{maintheorem}
\medskip

We remark that for large enough $\beta$ each $\sW_i$ has an overlap of length at least $1/2$ with each side-boundary and positive distance from the corners.
The fluctuations of the loop $\sL_0$ (the plateau at level $\lfloor H\rfloor$) from its limit $\sW_0$ along the side-boundaries are now addressed.

\medskip
\begin{maintheorem}
  [Cube-root Fluctuations] \label{mainthm-3}
  In the setting of Theorem~\ref{mainthm-2} suppose
 $\liminf_{n\to\infty}\alpha_n>\alpha_c$.
Consider $\sL_0^{(n)}$, the collection of macroscopic loops belonging to the level line at height $\lfloor H_n\rfloor$ of the SOS surface.
Then the maximal vertical fluctuation of $\sL_0^{(n)}$ from the boundary interval $[(\frac14L_n,0),(\frac34L_n,0)]$ is w.h.p.\ of order $L_n^{1/3+o(1)}$.
More precisely, let $\rho(x) = \min\{ y : (x,y)\in \sL_0^{(n)}\}$ be the vertical fluctuation of $\sL_0^{(n)}$ from the bottom boundary at coordinate $x$.
Then for any $\epsilon>0$ and large enough $n$,
\[ \pi_{\Lambda_n}\bigg( L_n^{1/3-\epsilon} < \sup_{\frac14 L_n \leq x \leq \frac34L_n} \rho(x)< L_n^{1/3+\epsilon} \bigg) > 1 - e^{-L_n^{\epsilon}}\,.\]
\end{maintheorem}
\medskip

As a direct corollary one can deduce an upper bound on the fluctuations of \emph{all} level lines.

\medskip
\begin{maincorollary}[Cascade of fluctuation exponents]\label{maincor-4}
In the same setting of Theorem~\ref{mainthm-3}, let $\rho(i,x)$ be the vertical fluctuation of $\sL_i^{(n)}$ from the bottom boundary at coordinate $x$.
Let $0<\xi<1$ and let $i = \lfloor \xi H_n \rfloor$. Then for any $\epsilon>0$,
\[ \lim_{n\to\infty} \pi_{\Lambda_n}\bigg( \sup_{\frac14 L_n \leq x \leq \frac34L_n} \rho(i,x) > L_n^{\frac{1-\xi}3+\epsilon} \bigg) = 0\,.\]
\end{maincorollary}
\medskip
\noindent\textbf{Proof.}
 Choose $\ell=(1+o(1))L^{1-\xi}$ such that $\lfloor H(\ell)\rfloor = \lfloor H_n \rfloor - i$ while $\alpha(\ell) > 2\alpha_c$ (recalling that $\alpha_c$ is arbitrarily small for large $\beta$).
Fix a coordinate $\frac14 L_n \leq x \leq \frac34L_n$ on the bottom boundary and consider the box $B=[x-\frac{\ell}2, x+\frac{\ell}2]\times [0, \ell]$. By imposing zero boundary conditions on the box $B$ we obtain a surface that is clearly stochastically dominated by the original one. By the choice of $\alpha(\ell)$, Theorem~\ref{mainthm-1} ensures w.h.p.\ the existence of a macroscopic loop $\sL_0'$ at height $\lfloor H(\ell)\rfloor$ for the new surface in $B$. Theorem~\ref{mainthm-3} guarantees that the vertical fluctuation of $\sL_0'$ at coordinate $x$ is at most ${\ell}^{1/3+\epsilon}$ except with probability $e^{-\ell^\epsilon}$. A union bound
extends this bound to the entire interval $x\in [\frac14 L_n, \frac34 L_n]$. Hence, w.h.p.\ the macroscopic loop $\sL_i$ of the original surface must be within $\ell^{1/3+\epsilon}=(1+o(1))L^{(1-\xi)/3+\epsilon}$ from the bottom interval.
\hfill\qed

Detailed proofs of Theorem~\ref{mainthm-1}--\ref{mainthm-3} will be given in a forthcoming paper.

\vspace{-0.75cm}
\section{Related work and open problems}
\vspace{-0.3cm}
\noindent \textbf{Limiting shape.}
In the two papers~\cite{ScSh1,ScSh2} the authors studied the limiting shape of the low temperature 2\Dim Ising with minus boundary under a prescribed small positive external field, proportional to the inverse of the side-length $L$. The behavior of the droplet of plus spins in this model is qualitatively similar to
the behavior of the top loop $\sL_0$ in our case. Here, instead of an external field, it is the subtle entropic repulsion phenomenon which induces the surface to rise to level $\lfloor H \rfloor$ and creates a macroscopic loop $\sL_0$.
In line with this connection, the description of the limiting shape of the plus droplet in the aforementioned works exactly coincides with our limit $\sW_0$ (with a different Wulff shape).

An important difference between the two models is of course the fact that in our case there exist $\lfloor H \rfloor \geq c \log L$ levels (rather than just one), which are interacting in two nontrivial ways. First, they are dependent as they are nested by definition. Second, they can weakly either attract or repel one another depending on the local geometry and height. Moreover, the box boundary itself can attract or repel the level lines. A prerequisite to proving Theorem~\ref{mainthm-2} is to overcome these ``pinning'' issues. As a side note we remark that at times such pinning issues have been overlooked in the relevant literature.

As for the fluctuations of the plus droplet from its limiting shape, it was argued in~\cite{ScSh1} that these should be normal (i.e., $\sqrt{L}$).
However, due to the analogy mentioned above between the models, it follows from our proof of Theorem~\ref{mainthm-3} that these fluctuations are in fact $L^{1/3+o(1)}$.

\smallskip
\noindent\textbf{Cube root fluctuations.}
It turns out that governing the behavior of the top level loop $\sL_0$
is a version of the usual measure associated with Ising-type contours (see, e.g.,~\cite{DKS}) crucially tilted by a factor of the form $\exp\big((\lambda/L)\mathrm{area}(\sL_0)\big)$ for
some fixed $\lambda>0$. Modulo this fact, whose proof is quite delicate, the $L^{1/3}$ fluctuations can be recovered by analyzing the behavior of the loop along mesoscopic boxes of dimensions $L^{2/3+o(1)} \times L^{1/3+o(1)}$, in the spirit of the approach of~\cite{DKS}.

There is a rich literature of contour models featuring similar distributions and cube root fluctuations. In some of these works (e.g.,~\cite{Alexander,FeSp,HV,Vel2}) the area term appears due to an externally imposed constraint (by conditioning and/or adding an external field). In others, modeling ordered random walks (e.g.,~\cite{CH,Johansson} to name a few), the area effect on a contour is due to the contours preceding it, and  the exact solvability of the model (e.g., via determinantal representations) plays an essential role in the analysis.

We stress that the area term driving the $L^{1/3}$ fluctuations in our setting arises naturally due to entropic repulsion.
By the lack of exact solvability for the $(2+1)\Dim$ SOS we must resort to cluster expansion techniques and contour analysis
as in the framework of~\cite{DKS}.
However, in the mesoscopic scale (required to establish the cube root fluctuations) this framework is particularly prone to the aforementioned pinning issues of level lines (to those above and below them, as well as the boundary), forming a major technical obstacle for the analysis.

\smallskip
\noindent\textbf{Open problems.} We conclude by mentioning two problems that remain unaddressed by our results. First, establish
the exponents for the fluctuations of all intermediate level lines from the side-boundaries.
We believe the upper bound in Corollary~\ref{maincor-4} features the correct cascade of exponents.
Second, find the correct fluctuation exponent of the level lines $\{\sL_i\}$ around the curved part of their limiting shapes $\{\sW_i\}$.

\end{document}